

\input amstex.tex
\input amsppt.sty

\TagsAsMath



\def\det{\operatorname{det}}

\def\ker{\operatorname{Ker}}
\def\supp{\operatorname{supp}}
\def\sup{\operatorname{sup}}
\def\inf{\operatorname{inf}}

\def\p{\partial}\def\at#1{\vert\sb{\sb{#1}}}
\def\pr{\sp{\prime}}

\def\half{\frac{1}{2}}

\def\vect#1{{\bold #1}}
\def\R{{\Bbb R}}  
\def\N{{\Bbb N}} \def\Z{{\Bbb Z}}

\def\Abs#1{\left\vert#1\right\vert}
\def\abs#1{\vert#1\vert}
\def\Norm#1{\left\Vert #1 \right\Vert}
\def\norm#1{\Vert #1 \Vert}
\def\tangent#1{T #1} \def\cotangent#1{T\sp{\ast} #1}

\def\ssb#1{\!\,\sb{\!\sb{#1}}}
\def\smsb#1{\sb{\sssize{#1}}}
\def\csb#1{\sb{\hskip -1pt #1}}
\def\cpr{\sp{\hskip -0.5pt {\sssize \prime}}}

\def\const{\operatorname{const}}

\def\sothat{\,\,\vert\,\,}
\def\habla#1{\nabla\sb{\hskip -2pt #1}}

\def\nut{\sb{\!o}}
\def\loc{\sb{\text{loc}}}
\def\comp{\sb{\text{comp}}}

\def\Month{\ifcase\month\or January\or February\or March\or
April\or May \or June \or
July\or August\or September \or October\or November\or
December\fi}

\def\authorinfo{%
\author
Andrew Comech
\medskip
{\rm Mathematics,
SUNY at Stony Brook,
${\Bbb {NY}}$  11794}
\endauthor
\email { comech${\sssize\@}$math.sunysb.edu}\endemail
}%

\def\proof{\demo{Proof} }
\def\endproof{\hfill$\square$\enddemo}
\def\sec#1{#1}

\def\ye#1{$\overset{\ssize #1}\to{\phantom{.}}$}
\def\noref#1\endref{}

\def\Without#1{}
\def\hvar{\hslash}
\def\thvar{\vartheta}
\def\sing#1{{S}\sb{#1}}

\magnification=1200
\hsize=4.6in
\vsize=7.0in
\hoffset=0.5in  
\voffset=0.0in
\nologo


\def\hvar{{\hslash}}

\def\Tbar{\overline{T}}

\document

\topmatter

\title
Optimal regularity\\
of Fourier integral operators\\
with one-sided folds
\endtitle

\authorinfo


\leftheadtext{Andrew Comech}
\rightheadtext{Fourier Integral Operators with one-sided folds }

\abstract
We obtain optimal continuity
in Sobolev spaces
for the Fourier integral operators
associated to singular canonical relations,
when one of the two projections 
is a Whitney fold.
The regularity depends on the type, $k$, of 
the other projection from the canonical relation
($k=1$ for a Whitney fold).
We prove that 
one loses $(4+\frac{2}{k})^{-1}$ 
of a derivative in the regularity properties.

The proof is based on the $L^2$ estimates
for oscillatory integral operators.

\endabstract


\endtopmatter 

\bigskip\head
\sec1. Introduction and results
\endhead

The Fourier integral
operators associated to singular canonical relations
(i.e., which are not local graphs) 
fall out of scope
of the classical theory of Fourier integral operators.
Their regularity properties are still to be studied.
The first step in this direction was due to the 
paper of R.B. Melrose and M.E. Taylor \cite{MeTa\ye{85}},
who showed that the canonical relation
with Whitney folds on both sides 
can be transformed to the canonical form.
(This idea originates
from the paper of Melrose \cite{Me\ye{76}}.)
Melrose and Taylor then derived 
the loss of $1/6$ of a derivative 
in the regularity properties 
of Fourier integral operators
with folding canonical relations,
vs. operators 
associated to local canonical graphs.

A corresponding result for oscillatory integral operators
(with not necessarily homogeneous phase functions)
was obtained by Y. Pan and C.D. Sogge \cite{PaSo\ye{90}},
who also relied on the reduction of the folding
canonical relation to the normal form.
An independent analytical approach 
to such operators in $\R^1$
was used by
D.H. Phong and E.M. Stein \cite{PhSt\ye{91}}.
This approach was generalized 
for operators in $\R^n$
by S. Cuccagna \cite{Cu\ye{97}},
who used fine almost orthogonal decompositions
of the integral kernels.
Let us also mention
thorough investigations
of oscillatory integral operators in $\R^1$
and related 
generalized Radon transforms in the plane
by Phong and Stein \cite{PhSt\ye{97}}
and by A. Seeger \cite{Se\ye{93}}, \cite{Se\ye{98}}.
We also refer to the survey of D.H. Phong \cite{Ph\ye{94}}.

The relation of the regularity properties 
of Fourier integral operators to the
rate of high-frequency decay of norms of
oscillatory integral operators 
was used
by A. Greenleaf and A. Seeger \cite{GrSe\ye{94}}
for deriving the {\it a priori} 
continuity of Fourier integral operators
associated to one-sided Whitney folds:
in general, one loses 
up to $1/4$ of a derivative
in the regularity properties.
Let us mention a recent result \cite{GrSe\ye{98}}
that for operators associated to canonical relations
with cusp singularities on one side there is a loss
of at most $1/3$ of a derivative.

In this paper, we are going to derive the optimal
regularity properties of the Fourier integral operators
associated to one-sided Whitney folds.
Let us recall the standard framework.
Let $X$ and $Y$ be $C^\infty$ manifolds 
of the same dimension $n$,
and let ${\Cal C}$ be a homogeneous canonical relation
${\Cal C}\subset\cotangent{X}\backslash 0
\times\cotangent{Y}\backslash 0$,
that is, 
${\Cal C}$ is lagrangian 
with respect to the 
difference of the canonical symplectic forms 
lifted from $\cotangent{X}$ and $\cotangent{Y}$
onto ${\Cal C}$.
If ${\Cal C}$ is locally a canonical graph,
that is, 
both projections
$$
\pi\ssb{L}:\;
{\Cal C}\to \cotangent{X}\backslash 0,
\qquad
\pi\ssb{R}:\;
{\Cal C}\to \cotangent{Y}\backslash 0
$$
are local diffeomorphisms,
then the regularity properties of 
Fourier integral operators 
associated to ${\Cal C}$
are well-known, cf. L. H\"ormander \cite{H\"o\ye{85}}:
Given a Fourier integral operator
${\goth F}\in I^m(X,\, Y,\,{\Cal C})$,
then for any real $s$
there is a continuous map
$
{\goth F}:\;H\comp^{s}(Y)
\to H\loc^{s-m}(X).
$
Here
$H^s(X)$
is the standard Sobolev space of order $s$.

Now let us state
the continuity of the Fourier integral operators
associated to 
canonical relations
with a Whitney fold on one side.
The continuity turns out to depend on 
the {\it type} 
of the projection
from the canonical relation
on the other side
(see the definition after the theorem).

\proclaim{Theorem \sec1.1}
Let
${\goth F}\in 
I^m(X,\, Y,\,{\Cal C})$ 
be a Fourier integral operator 
associated to the homogeneous
canonical relation 
${\Cal C}
\subset\cotangent{X}\backslash 0
\times\cotangent{Y}\backslash 0
$, 
such that one of the projections
${\Cal C}\to \cotangent{X}\backslash 0$,
$\;{\Cal C}\to \cotangent{Y}\backslash 0$
is a Whitney fold 
and the other is of type at most $k$ ($k=1$ for a Whitney fold).
Then, for any real $s$, ${\goth F}$ defines a continuous map
$$
{\goth F}:\;H\comp^{s}(Y)
\to
H\loc^{s-m-\,(4+\frac{2}{k})^{-1}}(X).
$$

This result is optimal,
in the sense that
there are operators associated
to the singular canonical relations
with one of the projections being a Whitney fold
and the other of type $k$,
when the smoothing stated above
can not be improved.
\endproclaim

\noindent
The above-mentioned case 
of two-sided Whitney folds 
corresponds to $k=1$.

The {\it type}
of a map of corank at most $1$ is to be defined as 
the highest order of vanishing
of the determinant of its Jacobian
in the directions of the kernel of its differential.
Let $M$ and $N$ be two $C^\infty$ manifolds 
of the same dimension
and let $\pi:\;M\to N$
be a smooth map.
We assume that 
the rank of $\pi$ drops simply by 1:
the corank of 
$d\pi$ is at most $1$ and 
the differential $d(\det d\pi)$ 
does not vanish in some neighborhood
of the critical variety
$$
\varSigma(\pi)
=\{p\in M\sothat \det d\pi\at{p}=0\}.
$$
Let 
$\vect{V}
\in C^\infty\left(\Gamma(\tangent{M}\at{U})\right)$ 
be a smooth vector field defined 
in some open neighborhood $U\subset M$ 
of a point $p\nut\in\varSigma(\pi)$,
which generates the kernel of $d\pi$:
$$
\vect{V}\at{U}\ne 0,
\qquad\qquad
\vect{V}\at{U\cap \varSigma(\pi)}\in\ker d\pi.
\tag{\sec1.1}
$$

\definition{Definition}
The type of $\pi$
at a point $p\nut\in\varSigma(\pi)$
is 
the smallest $k\in\N$ such that
$$
\vect{V}^k (\det d\pi)\at{p\nut}\ne 0.
\tag{\sec1.2}
$$
The type of $\pi$ at $p\in M\backslash\varSigma(\pi)$
is defined to be $0$.
\enddefinition

Since we assume that the rank
of $\pi$ drops simply and hence $d(\det d\pi)\ne 0$
on $\varSigma(\pi)$,
any other smooth vector field $\tilde\vect{V}$
which satisfies the conditions $(\sec1.1)$
can be represented 
in the neighborhood $U$
as
$\tilde\vect{V}=\varphi\cdot\vect{V}+(\det d\pi)\cdot\vect{W}$,
where 
$\varphi$ is a smooth function on $M$ 
which does not vanish on $\varSigma(\pi)$
and 
$\vect{W}$ is a smooth vector field.
(Let us note that $\det d\pi$ is only defined up to a
non-zero factor,
depending on the choice of local coordinates.)
As a consequence,
the above Definition does not
depend on the choice of $\vect{V}$.

\demo{Remark}
In the context of singular integral operators,
the type conditions 
can be defined without the assumption
that the rank drops simply,
see for example
\cite{PhSt\ye{91}},
\cite{Se\ye{93}},
\cite{PhSt\ye{94}},
\cite{Co\ye{97}},
and \cite{Se\ye{98}}.
We use the above Definition
since
the assumption
of Theorem \sec1.1 
that one of the projections from 
the canonical relation
is a Whitney fold guarantees that
the other projection drops its rank simply by 1:
This is because
the corank of both projections is the same
(and hence not greater than $1$)
and 
$\det d\pi\ssb{L}$, $\det d\pi\ssb{R}$
are equal up to a non-zero factor
(and both vanish simply on the critical 
variety).

\enddemo

An example of a map of type at most $k$
is a map with a Morin $\sing{1\csb{k}}$-singularity,
with $k$ units
(see R. Thom \cite{Th\ye{63}} 
and B. Morin \cite{Mo\ye{65}}).
For example, 
type $k=1$ 
unambiguously corresponds to the Whitney fold.
$k=2$ for the {\it Whitney Pleat}, 
or the {\it Simple Cusp},
at the cusp point
(and $k=1$ at neighboring critical points);
$k=3$ for the {\it Swallow Tail}
at the ``swallow tail point'',
et cetera.

\subhead{Relation with oscillatory integral operators}
\endsubhead
We will use the results of Greenleaf and Seeger
on the relation of Fourier integral operators
and oscillatory integral operators.
As they showed in \cite{GrSe\ye{94}},
the regularity in Sobolev spaces of 
{\it singular Fourier integral operators} 
associated to singular canonical relations
is determined by the rate
of decay of the $L^2$-operator norm 
of {\it oscillatory integral operators}
associated to similar canonical relations.

The oscillatory integral operators are of the form
$$
T\sb{\lambda}u(x)=\underset{\R^n}\to\int
e^{i\lambda S(x,\thvar)}
\,\psi(x,\thvar)\,u(\thvar)\, d\thvar,
\quad
\psi\in C\comp^\infty(\R^n\times\R^n),
\ \lambda\gg 1.
\tag{\sec1.3}
$$
We will write subscripts 
for two copies of $\R^n$:
$x\in\R^n\smsb{L}$, $\thvar\in\R^n\smsb{R}$.
The canonical relation 
${\Cal C}$ associated to 
the oscillatory integral operator $(\sec1.3)$
is given by
$$
{\Cal C}=\{
(x,S\sb{x})\times(\thvar,S\sb{\thvar})\sothat
x\in\R^n\smsb{L},\,\thvar\in\R^n\smsb{R}\}
\subset\cotangent\R^n\smsb{L}\!\times\cotangent\R^n\smsb{R}.
$$
Here $S\sb{x}$ stands for the components of 
the 1-form
$d\sb{x}S$, etc.
Using the isomorphism
$
\R^n\smsb{L}\!\times\R^n\smsb{R}\cong
{\Cal C},
$
we write the projections from the canonical relation
onto the first and second factors 
of $\cotangent\R^n\smsb{L}\!\times\cotangent\R^n\smsb{R}$
in the following form:
$$
\pi\ssb{L}:\;(x,\,\thvar)\mapsto(x,\,S\sb{x}),
\qquad\qquad
\pi\ssb{R}:\;(x,\,\thvar)\mapsto(\thvar,\,S\sb{\thvar}).
\tag{\sec1.4}
$$
The projections $(\sec1.4)$ 
are degenerate
on the variety
where the determinant of the mixed Hessian of $S$ vanishes.
We will use the notation
$$
h(x,\thvar)=\det S\sb{x\thvar},
$$
so that the common critical variety
of the projections $\pi\ssb{L}$, $\pi\ssb{R}$
is given by
$$
\varSigma=\left\{\,(x,\thvar) 
\sothat
h(x,\thvar)=0 \right\}.
$$

\proclaim{Theorem \sec1.2}
Let $T\sb{\lambda}$ be an oscillatory integral operator 
$(\sec1.3)$ 
with a smooth compactly supported density $\psi(x,\thvar)$
and a smooth
(not necessarily homogeneous)
phase function 
$S(x,\thvar)$, $x$, $\thvar\in\R^n$,
such that one of the projections
$\pi\ssb{L}:\;(x,\,\thvar)\mapsto (x,\,S\sb{x})$
and 
$\pi\ssb{R}:\;(x,\,\thvar)
\mapsto (\thvar,\,S\sb{\thvar})$
from the associated canonical relation
is a Whitney fold and the other 
is of type at most $k$.

Then there is the following estimate on 
$T\sb{\lambda}$:
$$
\Norm{T\sb{\lambda}}\sb{L^2\to L^2}
\le\const\lambda^{-\frac{n}{2}+(4+\frac{2}{k})^{-1}}.
\tag{\sec1.5}
$$
This result is optimal
(see Section \sec4).
\endproclaim

Let us mention here
that the most general estimates
for oscillatory integral operators
with real analytic phase functions
in $\R^1$
are derived in \cite{PhSt\ye{97}}.

According to Greenleaf and Seeger \cite{GrSe\ye{94}},
Theorem \sec1.1
follows from Theorem \sec1.2.
The proof of Theorem \sec1.2 
is in Sections \sec2
and \sec3.
The sharpness of the results 
is discussed in Section \sec4.

\bigskip\head
\sec2. Dyadic creeping to the critical variety
\endhead

We use the dyadic 
decomposition
$
\sum\sb{N\in\Z}
\beta(2^N t)=1
$
(for $t>0$),
where $\beta(t)$ 
is a smooth function
supported in $[1/2,2]$,
$0\le\beta(t)\le 1$, $\beta(t)\equiv 1$
in a neighborhood of $t=1$.
We put $\beta\sb{+}(t)\equiv\beta(t)$,
$\beta\sb{-}(t)\equiv\beta(-t)$,
to take care of positive and negative values separately
(usually we will not write $\pm$-subscripts).
We define $\bar{\beta}\in C\comp^\infty([-2,2])$
by 
$\bar{\beta}(t)=\beta(\abs{t})$ for $\abs{t}\ge 1$,
$\bar{\beta}(t)\equiv 1$ for $\abs{t}\le 1$.
There is the following 
partition of 1 which depends on $N\nut\in\Z$:
$$
1
=\sum\sb{\pm}\sum\sb{N\in \Z,\,N<N\nut}
{\beta\sb{\pm}}(2^{N}h(x,\thvar))
+{\bar{\beta}}(2^{N\nut}h(x,\thvar)).
\tag{\sec2.1}
$$

We define the operators
$T\sb{\lambda}\sp{\pm\hvar}$ and $\Tbar\sb{\lambda}\sp{\hvar}$ by
$$
\align
T\sb{\lambda}\sp{\pm\hvar} u(x)
&=\int\sb{\R\ssb{R}^n}
\,e^{i\lambda S(x,\thvar)} 
\,\psi(x,\thvar)
\,\beta\sb{\pm}(\hvar^{-1}h(x,\thvar))
\,u(\thvar)\,d\thvar,
\tag{\sec2.2}
\\ \\
{
\Tbar\sb{\lambda}\sp{\hvar}} u(x)
&=\int\sb{\R\ssb{R}^n}
\,e^{i\lambda S(x,\thvar)} 
\,\psi(x,\thvar)\,
{\bar{\beta}}(\hvar^{-1}h(x,\thvar))
\,u(\thvar)\,d\thvar,
\tag{\sec2.3}
\endalign
$$
(here $\hvar=2^{-N}$ refers to the magnitude
of $h(x,\thvar)$), 
and decompose $(\sec1.3)$ into a sum
$$
T\sb{\lambda}
=
\sum\sb{\pm}
\sum\sb{\hvar>\hvar\nut(\lambda)}\sp{\hvar\le 2D}
T\sb{\lambda}\sp{\pm\hvar}
+{\Tbar\sb{\lambda}\sp{\hvar\nut(\lambda)}},
\quad\hvar=2^{-N}, 
\quad
N\in\Z,
\tag{\sec2.4}
$$
where 
$D$ is the uniform bound on 
$\abs{h(x,\thvar)}$,
and $\hvar\nut(\lambda)$
will be chosen 
so that the optimal estimates on 
$T\sb{\lambda}\sp{\pm\hvar}$
and 
$\Tbar\sb{\lambda}\sp{\hvar}$
coincide
when $\hvar=\hvar\nut(\lambda)$.
It suffices to consider the
estimates on 
$T\sb{\lambda}\sp{\pm\hvar}$,
$\Tbar\sb{\lambda}\sp{\hvar}$
for $\hvar< 1$.

\bigskip

We will 
prove the estimates
$\norm{T\sb{\lambda}\sp{\pm\hvar}}\sb{L^2\to L^2}
\le\const\lambda^{-\frac{n}{2}}\hvar^{-\half}$
(Theorem \sec2.1)
and
$\norm{\Tbar\sb{\lambda}\sp{\hvar}}\sb{L^2\to L^2}
\le\const\lambda^{-\frac{n-1}{2}}\hvar^{\half+\frac{1}{2k}}$
(Theorem \sec3.1);
these estimates meet at
$$
\hvar\nut(\lambda)=\lambda^{-\frac{k}{2k+1}}.
\tag{\sec2.5}
$$
Using the corresponding estimates for
the operators in the right-hand side of $(\sec2.4)$,
we arrive at the statement of Theorem \sec1.2.

\proclaim{Theorem \sec2.1}
Let
the projection 
$\pi\ssb{L}:\;(x,\thvar)\mapsto(x,S\sb{x})$ 
be a Whitney fold
and let the projection 
$\pi\ssb{R}:\;(x,\thvar)\mapsto(\thvar,S\sb{\thvar})$ 
be of type at most $k$,
for some $k\in\N$.
Then, as long as $\hvar\ge\lambda^{-\half}$,
there is the following estimate:
$$
\Norm{T\sb{\lambda}\sp{\pm\hvar}}\sb{L^2\to L^2}
\le \const\lambda^{-\frac{n}{2}}\hvar^{-\half}.
\tag{\sec2.6}
$$
\endproclaim

Since both $T\sb{\lambda}\sp{\hvar}$ 
and $T\sb{\lambda}\sp{-\hvar}$ require the
same argument, we will always restrict
the consideration to $T\sb{\lambda}\sp{\hvar}$.
Also, unless otherwise stated, 
the norm $\norm{\ }$ will refer to the 
$L^2$ operator norm.

\demo{Remark}
According to $(\sec2.5)$, we are only interested
in
$\hvar\ge\lambda^{-\frac{k}{2k+1}}$;
we will only give the proof for this case.
This proof 
has already appeared 
in the author's paper \cite{Co\ye{97}},
but 
we reproduce it for the sake of completeness.

The proof for 
$\hvar\ge\lambda^{-\half}$
can be obtained by some elaboration
of almost orthogonal decompositions.
The estimate $(\sec2.6)$ 
becomes useless for $\hvar<\lambda^{-\half}$.
\enddemo

The proof involves the spatial 
decomposition with respect to $\thvar$,
with the step $\hvar$.
We use the notation 
$\left(T\sb{\lambda}\sp{\hvar}\right)\sb{\varTheta}$
for $T\sb{\lambda}\sp{\hvar}$ localized near the point 
$\hvar\varTheta\in\R^n\smsb{R}$;
here $\varTheta$ is a point on the integer lattice $\Z^n$.

As long as $\hvar\ge\lambda^{-\half+\varepsilon}$,
$\varepsilon>0$,
the argument similar to the one used by
S. Cuccagna \cite{Cu\ye{97}} shows that
$\pi\ssb{L}$ being a Whitney fold is a sufficient condition
for the pieces 
$\left(T\sb{\lambda}\sp{\hvar}\right)\sb{\varTheta}$
to be almost orthogonal with respect to different 
values of $\varTheta\in\Z^n$:
$$
\Norm{
\left(T\sb{\lambda}\sp{\hvar}\right)\sb{\varTheta}
\left(T\sb{\lambda}\sp{\hvar}\right)\sb{W}\sp{\ast}
},\,\,
\Norm{
\left(T\sb{\lambda}\sp{\hvar}\right)\sb{\varTheta}\sp{\ast}
\left(T\sb{\lambda}\sp{\hvar}\right)\sb{W}
}
\le\const\lambda^{-n}\hvar^{-1}\Abs{\varTheta-W}^{-N},
$$
for any $N\in\N$.
Here $\Abs{\varTheta-W}$ is the distance
between the points $\varTheta$, $W$ in $\Z^n$.
Then, the Cotlar-Stein almost orthogonality
lemma \cite{St\ye{93}} applies.

Let us derive
the individual estimates on 
$\left(T\sb{\lambda}\sp{\hvar}\right)\sb{\varTheta}$,
which are similar to H\"ormander's estimates
for non-degenerate oscillatory integral operators
in $\R^n$.
We consider the integral kernel 
of the composition
$\left(T\sb{\lambda}\sp{\hvar}\right)\sb{\varTheta}
\left(T\sb{\lambda}\sp{\hvar}\right)\sb{\varTheta}\sp{\ast}$:
$$
K\left(
\left(T\sb{\lambda}\sp{\hvar}\right)\sb{\varTheta}
\left(T\sb{\lambda}\sp{\hvar}\right)\sb{\varTheta}\sp{\ast}
\right)(x,y)
=\int d^n\thvar\,
e^{i\lambda(S(x,\thvar)-S(y,\thvar))}
\chi(\hvar^{-1}\thvar-\varTheta)
\times\dots\, .
\tag{\sec2.7}
$$
We integrate by parts in $(\sec2.7)$, using the 
operator
$$
L\sb{\thvar}=\frac{1}{i\lambda}
\cdot
\frac{(S\sb{\thvar}(x,\thvar)-S\sb{\thvar}(y,\thvar))
\cdot\habla{\thvar} }
{\Abs{S\sb{\thvar}(x,\thvar)-S\sb{\thvar}(y,\thvar) }^2 }.
$$
When acting on cut-offs,
$\habla{\thvar}$ contributes $\hvar^{-1}$
(we will discuss this below in more details).
This is the pay for decomposing $T\sb{\lambda}$ 
with respect to the values of $h(x,\thvar)$; 
one takes over when integrating
with respect to $\thvar$ in $(\sec2.7)$,
since $\left(T\sb{\lambda}\sp{\hvar}\right)\sb{\varTheta}$
has the support of size $\hvar$ in $\thvar$-variables. 

We then apply the Schur lemma,
i.e., integrate with respect to $x$ (or $y$)
the absolute value of $(\sec2.7)$ with the extra factor 
$\left(1+\lambda\hvar
\Abs{S\sb{\thvar}(x,\thvar)-S\sb{\thvar}(y,\thvar)}
\right)^{-N}$.
It is convenient to change the variables of integration:
$x\mapsto \eta=S\sb{\thvar}(x,\thvar)$,
$\dsize
\int \!\!dx \,d\thvar\to
\int\!\!
\frac{d \eta\, d\thvar}{\abs{\det S\sb{x\thvar}}}.
$
The integration with respect to $\eta$
contributes $(\lambda\hvar)^{-n}$,
the
integration with respect to $\thvar$ 
contributes $\hvar^{n}$,
and $\abs{\det S\sb{x\thvar}}\approx\hvar$.
This yields the estimate
$\norm{
\left(T\sb{\lambda}\sp{\hvar}\right)\sb{\varTheta}
\left(T\sb{\lambda}\sp{\hvar}\right)\sb{\varTheta}^\ast
}\le\const\lambda^{-n}\hvar^{-1},
$
and hence proves the bound $(\sec2.6)$
on 
$\left(T\sb{\lambda}\sp{\hvar}\right)\sb{\varTheta}$.

We need to control that,
during the integrations by parts,
the derivative $\habla{\thvar}$ 
contributes at most $\const\hvar^{-1}$ even when it acts 
on the denominator of $L\sb{\thvar}$ itself.
For this, the map 
$\pi\ssb{R}\at{\thvar}:\,x\mapsto S\sb{\thvar}(x,\thvar)$
needs to satisfy certain convexity condition
on the support of
$T\sb{\lambda}\sp{\hvar}$:
$$
\Abs{S\sb{\thvar}(x,\thvar)-S\sb{\thvar}(y,\thvar)}
\ge\const\hvar\,\abs{x-y}.
\tag{\sec2.8}
$$
As the matter of fact, $\pi\ssb{R}$
does not generally satisfy the condition $(\sec2.8)$
on the entire support of $T\sb{\lambda}\sp{\hvar}$,
and we are going to introduce one more localization.

In the rest of this section, we discuss the
construction of this new localization
(which will also be used in Section \sec3)
and prove that $(\sec2.8)$ is valid on the support
of each of the pieces of $T\sb{\lambda}\sp{\hvar}$.

Since one of the projections 
from the associated canonical relation
is a Whitney fold
(we assume it is $\pi\ssb{L}$),
the rank of the mixed Hessian
$S\sb{x\thvar}$ 
on the critical variety
is equal to $n-1$.
We choose local coordinates $x=(x\pr,x\sb{n})$
and $\thvar=(\thvar\pr,\thvar_n)$
so that $S\sb{x\cpr\thvar\cpr}$
is non-degenerate.
We introduce 
the vector field 
$\vect{K}\ssb{R}$,
$$
\vect{K}\ssb{R}
=\p\sb{x\csb{n}}
-
S\sb{x\sb{n}\thvar\cpr}(x,\thvar)
S\sp{\thvar\cpr x\cpr}(x,\thvar)
\p\sb{x\cpr};
\tag{\sec2.9}
$$
we wrote
$S\sp{\thvar\cpr x\cpr}(x,\thvar)$
for the inverse to the matrix
$S\sb{x\cpr \thvar\cpr}$ at a point $(x,\thvar)$.
It can be checked immediately that 
this vector field satisfies
$$
\vect{K}\ssb{R}\at{\varSigma}
\in\ker d\pi\ssb{R}.
$$
Since the type of $\pi\ssb{R}$ is not greater
than $k$,
we can assume,
in the agreement with the definition
of type of the map (see $(\sec1.2)$),
that on the support of the integral
kernel of $T\sb{\lambda}$
$$
\abs{\vect{K}\ssb{R}^{k\cpr}h(x,\thvar)}\ge \kappa>0,
\tag{\sec2.10}
$$
for some positive constant $\kappa$
and for some integer $k\pr\le k$.
For definiteness,
we assume that $k\pr=k$
(this is the ``worst'' case).

We fix two smooth functions
$\rho\sb{-}$ and $\rho\sb{+}$,
supported in $(-\infty,1]$ and $[-1,\infty)$,
respectively,
such that
$
\rho\sb{-}(t)+\rho\sb{+}(t)=1,
$
$\,t\in\R,$
and define the following partition of 1:
$$
1=\sum\sb{\sigma}\rho\sp{\hvar}\sb{\sigma}(x,\thvar),
\qquad
\sigma=(\sigma\sb{1},\,\dots,\,\sigma\sb{k-1}),
\quad
\sigma\sb{j}=\pm 1,
$$
where
$$
\rho\sp{\hvar}\sb{\sigma}(x,\thvar)
\equiv\prod\sb{j=1}^{k-1}
\rho\sb{\sigma\sb{j}}(\hvar^{-1}
\vect{K}\ssb{R}^{j}h(x,\thvar)).
\tag{\sec2.11}
$$
We then split $T\sb{\lambda}\sp{\hvar}$
into $\sum\sb{\sigma}T\sb{\lambda,\sigma}\sp{\hvar}$,
multiplying the integral kernel of $T\sb{\lambda}\sp{\hvar}$
by the functions 
$\rho\sp{\hvar}\sb{\sigma}(x,\thvar)$:
$$
T\sb{\lambda,\sigma}\sp{\hvar}u(x)
=\int e^{i\lambda S(x,\thvar)}
\beta(\hvar^{-1}h)\rho\sp{\hvar}\sb{\sigma}(x,\thvar)
\psi(x,\thvar) u(\thvar) d\thvar.
\tag{\sec2.12}
$$

\proclaim{Proposition \sec2.2}
The map 
$\pi\ssb{R}\at{\thvar}:\,x\mapsto S\sb{\thvar}(x,\thvar)$
satisfies the convexity condition
$(\sec2.8)$
on the support of the integral kernel of each 
$T\sb{\lambda,\sigma}\sp{\hvar}$.

\endproclaim

This proposition allows the integration 
by parts in $(\sec2.7)$,
and thus finishes the proof of Theorem \sec2.1.
The proof of the proposition itself is based 
on two lemmas below.

We consider the map 
$\pi\ssb{R}\at{\thvar}:\;
x\mapsto \eta=S\sb{\thvar}(x,\thvar)$
as the composition
$$
\pi\ssb{R}\at{\thvar}:\;
x\overset{\pi\sp{\prime}}\to\longmapsto
(\eta\pr=S\sb{\thvar\cpr},x\sb{n})
\overset{\pi\sp{s}}\to\longmapsto
(\eta\pr,\eta_n=S\sb{\thvar\csb{n}}).
\tag{\sec2.13}
$$
Here
$\eta_n$ is considered
as a function of $\eta\pr$ and $x\sb{n}$: 
$
\eta\sb{n}(S\sb{\thvar\cpr}(x,\thvar),x\sb{n})
=S\sb{\thvar\csb{n}}(x,\thvar).
$

\demo{Remark}
The kernel of the differential $d\pi\sp{s}$
is certainly generated by the vector
$\left(\p\sb{x\csb{n}}\right)\sb{\eta\cpr}$
(the subscript refers to choosing 
$\eta\pr$ and $x\sb{n}$
as the independent variables),
and there is a convenient relation
$$
\vect{K}\ssb{R}
=\left(\p\sb{x\csb{n}}\right)\sb{\eta\pr},
$$
which motivated the definition $(\sec2.9)$ 
of $\vect{K}\ssb{R}$.
\enddemo

Since $\det S\sb{x\cpr\thvar\cpr}\ne 0$,
the map $\pi\sp{\prime}$ in $(\sec2.13)$ 
is a diffeomorphism (at least locally),
and hence we may assume that it satisfies 
$$
\abs{
\pi\pr(x)
-\pi\pr(y)}\ge\const\abs{x-y}.
\tag{\sec2.14}
$$

Now we work in the $(\eta\pr,x\sb{n})$-space;
we need to show that
$\pi\sp{s}$ in $(\sec2.13)$
satisfies
$$
\abs{
\pi\sp{s}(\eta\pr,x\sb{n})
-\pi\sp{s}(\zeta\pr,y_n)}
\ge\const\hvar\cdot\text{dist}[(\eta\pr,x\sb{n}),(\zeta\pr,y_n)],
\tag{\sec2.15}
$$
for appropriate ranges of the values of 
$\eta\pr$, $\zeta\pr$, $x\sb{n}$, and $y_n$.

We denote by ${\Cal L}$
the line segment 
from the point $(\eta\pr,x\sb{n})$
to $(\zeta\pr,y_n)$.
Since the first $n-1$ components of $\pi\sp{s}$
are identities, the inequality $(\sec2.15)$ 
is trivially satisfied if ${\Cal L}$
is outside the conic neighborhood 
of magnitude $c\hvar$
(where $c>0$ is to be chosen later)
of the directions $\pm(\p\sb{x\csb{n}})\sb{\eta\pr}$
in the $(\eta\pr,x\sb{n})$-space.

Now let ${\Cal L}$ be inside the $c\hvar$-cone
around $\pm(\p\sb{x\csb{n}})\sb{\eta\pr}$; then
the value of 
$\abs{\eta\pr-\zeta\pr}$ is bounded by
$c\hvar\abs{x\sb{n}-y_n}$.
According to the Mean Value theorem 
applied to $\eta_n(\eta\pr,x\sb{n})$,
there is the following bound from below
for the left-hand side of $(\sec2.15)$:
$$
\aligned
&\abs{\eta\sb{n}(\eta\pr,x_n)-\eta\sb{n}(\zeta\pr,y_n)}
\\
&\ge
\abs{x\sb{n}-y\sb{n}}
\cdot
\underset{{\Cal L}}\to\inf
\abs{\left(\p\sb{x\csb{n}}\right)\sb{\eta\pr}
\eta\sb{n}}
-
\abs{\eta\pr-\zeta\pr}\cdot
\underset{{\Cal L}}\to\sup
\abs{\nabla\sb{\eta\pr}\eta\sb{n}}
\\
&\ge
\abs{x\sb{n}-y\sb{n}}
\cdot
\left(
\underset{{\Cal L}}\to\inf
\abs{
 \left(\p\sb{x\csb{n}}\right)\sb{\eta\pr}\eta\sb{n}
}
-
c\hvar\,
\underset{{\Cal L}}\to\sup
\abs{\nabla\sb{\eta\pr}\eta\sb{n}}
\right).
\endaligned
\tag{\sec2.16}
$$
If we show 
that 
$
\underset{{\Cal L}}\to\inf
\abs{
 \left(\p\sb{x\csb{n}}\right)\sb{\eta\pr}\eta\sb{n}
}
$
is of magnitude $\hvar$,
then we may choose $c$ sufficiently small
so that
the inequality $(\sec2.8)$ follows.

The value of the derivative
$\left(\p\sb{x\csb{n}}\right)\sb{\eta\pr}\!\eta\sb{n}$
can be determined from 
the decomposition
$\pi\ssb{R}\at{\thvar}
=\pi\sp{s}\circ\pi\sp{\prime}$.
Considering the determinants 
of the Jacobi matrices,
$
J(\pi\ssb{R}\at{\thvar})=J(\pi\sp{s})\cdot J(\pi\sp{\prime}),
$
we obtain
$
h(x,\thvar)=
\left(\p\sb{x\csb{n}}\right)\csb{\eta\cpr}\!\eta\sb{n}
\cdot\det S\sb{x\cpr\thvar\cpr}.
$
Hence,

\proclaim{Lemma 1}
There is the relation
$
\left(\p\sb{x\csb{n}}\right)\sb{\eta\pr}
\!\eta\sb{n}
=\frac{h(x,\thvar)}{\det S\sb{x\cpr\thvar\cpr}}.
$
\endproclaim

Now we only need to check that
$\ h\ge\const\hvar$
everywhere on ${\Cal L}$,
if the length of ${\Cal L}$,
$\abs{{\Cal L}}\equiv
\operatorname{dist}[(\eta\pr,x\sb{n}),(\zeta\pr,y_n)]$,
is sufficiently small.

\proclaim{Lemma 2}
If $\abs{{\Cal L}}\le \frac{1}{12}$,
then $h\ge\frac{\hvar}{4}$
everywhere on ${\Cal L}$.

\endproclaim

\noindent
We thus admit that the line segment
${\Cal L}$
could be not entirely on the support 
of the integral kernel of $T\sb{\lambda}\sp{\hvar}$,
where $h\ge\hvar/2$.

\proof
Since both $(\eta\pr,x\sb{n})$ and $(\eta\pr,y_n)$
are on the support of 
the integral kernel of 
the operator
$T\sb{\lambda,\sigma}\sp{\hvar}$
defined by $(\sec2.12)$,
we have
$$
h\ge\frac{\hvar}{2}
\qquad\text{at the points $(\eta\pr,x\sb{n})$ and $(\zeta\pr,y_n)$},
\tag{\sec2.17}
$$
$$
\sigma\sb{j}
\cdot \vect{K}\ssb{R}^{j}h 
\ge -\hvar,
\quad j<k,
\qquad\text{at the points $(\eta\pr,x\sb{n})$ and $(\zeta\pr,y_n)$},
\tag{\sec2.18}
$$
for all $j<k$.
Also, 
according to $(\sec2.10)$,
$$
\abs{\vect{K}\ssb{R}^{k} h}\ge\kappa
\qquad\text{everywhere.}
\tag{\sec2.19}
$$

Let $t$ be a parameter on the line segment ${\Cal L}$,
changing from $t=0$ at the point $(\eta\pr,x\sb{n})$
to 
$t=\abs{\Cal L}$
at the point 
$(\zeta\pr,y_n)$;
$\p_t=\left(\p\sb{x\csb{n}}\right)\sb{\eta\pr}$.
We consider $h\at{\Cal L}$ 
as a function of $t$.
As long as
${\Cal L}$
is in the $c\hvar$-cone around $\pm (\p\sb{x\csb{n}})\sb{\eta\pr}$,
$$
\p\sb{t}^{j}h
=\vect{K}\ssb{R}\sp{j}h(x,\thvar)\at{\Cal L}
\qquad
\text{modulo terms of magnitude \ $c\hvar$,}
\tag{\sec2.20}
$$
for any $j\le k$.

If $c$ is sufficiently small,
then 
due to the inequalities $(\sec2.17)$-$(\sec2.20)$
we have
$$
h(0)\ge\frac{\hvar}{2},
\qquad
h(\abs{{\Cal L}})\ge\frac{\hvar}{2},
$$
$$
\sigma\sb{j}h^{(j)}(0)\ge -2\hvar,
\qquad
\sigma\sb{j}h^{(j)}(\abs{{\Cal L}})\ge -2\hvar,
\qquad\text{for all $j<k$},
\tag{\sec2.21}
$$
and also 
$$
\sigma_k h^{(k)}(t)>-3\hvar,
\qquad\text{for all $0\le t\le \abs{{\Cal L}}$,}
\tag{\sec2.22}
$$
where $\sigma_k$ is equal to $1$ or $-1$.
For our convenience, we have weakened the bound
in the right-hand side of $(\sec2.22)$.
We will base the rest of the argument
on the following elementary inequality:

\proclaim{Lemma}
Let $f(t)\in C^{1}([0,l])$.
If there is a uniform bound
$
\ \sigma f\pr(t)\ge -\epsilon,
$
where $\sigma$ is a constant equal to $\pm 1$ 
and $\epsilon>0$, then
$$
\operatorname{min}[f(0),f(l)]-\epsilon l
\le f(t)\le
\operatorname{max}[f(0),f(l)]+\epsilon l,
\qquad
\text{for any $0\le t\le l$.}
$$

\endproclaim

  From $(\sec2.22)$ 
and from the above Lemma 
(where we take $\epsilon=3\hvar$)
we conclude that
$\sigma\sb{k-1}h\sp{(k-1)}(t)
\ge -2\hvar-2\hvar\abs{{\Cal L}}
\ge -3\hvar$,
for any $t$ between $0$ and $\abs{{\Cal L}}$.
Continuing by induction, we conclude that 
$\sigma\sb{1}h\pr(t)\ge -3\hvar$.
Therefore, again from the above Lemma,
we deduce that
everywhere between $0$ and $\abs{{\Cal L}}$
$h(t)\ge\frac{\hvar}{2}-3\hvar\abs{{\Cal L}}$,
and this is not less than $\frac{\hvar}{4}$
as long as $\abs{{\Cal L}}\le\frac{1}{12}$.
\endproof

\bigskip\head
\sec3. Almost orthogonal decompositions 
near the critical variety
\endhead

Now we consider the operator $\Tbar\sb{\lambda}\sp{\hvar}$
defined by $(\sec2.3)$:
$$
{\Tbar\sb{\lambda}\sp{\hvar}} u(x)
=\int
\,e^{i\lambda S(x,\thvar)} 
\,\psi(x,\thvar)\,
{\bar{\beta}}(\hvar^{-1}h(x,\thvar))
\,u(\thvar)\,d\thvar,
\quad
\psi\in C\comp^\infty(\R^n\smsb{L}\times\R^n\smsb{R}),
$$
where
$\bar{\beta}\in C\comp^\infty(\R)$,
$\supp\bar{\beta}\subset[-2,2].$
The support of this operator contains the critical variety
$\varSigma=\{ \det S\sb{x\thvar}(x,\thvar)=0\}$.

\proclaim{Theorem \sec3.1}
If the projection $\pi\ssb{L}$ is a Whitney fold
and the projection $\pi\ssb{R}$ 
is of type at most $k$,
then,
as long as 
$
\hvar\ge\lambda^{-\frac{1}{2}},
$
there is the following estimate:
$$
\norm{\Tbar\sb{\lambda}\sp{\hvar}}\sb{L^2\to L^2}
\le \const
\lambda^{-\frac{n-1}{2}}\hvar^{\frac{1}{2}+\frac{1}{2k}}.
\tag{\sec3.1}
$$

\endproclaim

\demo{Remark}
According to $(\sec2.5)$, 
we only need the estimate $(\sec3.1)$
for 
$\hvar=\lambda^{-\frac{k}{2k+1}}$.
We will prove
Theorem \sec3.1
assuming that
$$
\hvar\ge\lambda^{-\frac{k}{2k+1}},
\tag{\sec3.2}
$$
to avoid unnecessary details.
\enddemo

We start with the decomposition
of
$\Tbar\sb{\lambda}\sp{\hvar}$
into
$$
\Tbar\sb{\lambda}\sp{\hvar}
=\sum\sb{\sigma}\Tbar\sb{\lambda,\sigma}\sp{\hvar},
\qquad\sigma=(\sigma_1,\,\dots,\,\sigma\sb{k-1}),
\quad \sigma_j=\pm 1,
$$
with respect to the signs of the derivatives
$\vect{K}\ssb{R}\sp{j}h$, $1\le j\le k-1$:
We use introduced earlier functions
$\rho\sp{\hvar}\sb{\sigma}(x,\thvar)$
(see $(\sec2.11)$) and define
$$
\Tbar\sb{\lambda,\sigma}\sp{\hvar}u(x)
=\int e^{i\lambda S(x,\thvar)}
\bar\beta(\hvar^{-1}h)\rho\sp{\hvar}\sb{\sigma}(x,\thvar)
\psi(x,\thvar) u(\thvar) d\thvar.
\tag{\sec3.3}
$$
On the support of the integral kernel of
$\Tbar\sb{\lambda,\sigma}\sp{\hvar}$
the following inequalities
are satisfied:
$$
\sigma\sb{j}\vect{K}\ssb{R}\sp{j}h(x,\thvar)
\ge -\hvar,
\tag{\sec3.4}
$$
for all $j$ between $1$ and $k-1$.
Let us mention that for $k=1$ no decomposition is needed.

We will consider the operators
$\Tbar\sb{\lambda,\sigma}\sp{\hvar}$ 
with different sets $\sigma$
separately.
Given $\sigma$,
we decompose the corresponding
$\Tbar\sb{\lambda,\sigma}\sp{\hvar}$
into
$$
\Tbar\sb{\lambda,\sigma}\sp{\hvar}
=\sum\sb{X\in\Z^n}\sum\sb{\varTheta\in\Z^n}
\left(\Tbar\sb{\lambda,\sigma}\sp{\hvar}
\right)\sb{X\varTheta},
\tag{\sec3.5}
$$
where
$
\left(\Tbar\sb{\lambda,\sigma}\sp{\hvar}
\right)\sb{X\varTheta}
$
is an operator with the integral kernel
$$
\chi(\hvar^{-\frac{1}{k}}x-X)
\cdot
K\left(\Tbar\sb{\lambda,\sigma}\sp{\hvar}
\right)(x,\thvar)
\cdot\chi(\hvar^{-1}\thvar-\varTheta).
\tag{\sec3.6}
$$
Here 
$K\left(\Tbar\sb{\lambda,\sigma}\sp{\hvar}
\right)(x,\thvar)$
stands for the integral kernel of
$\Tbar\sb{\lambda,\sigma}\sp{\hvar}$.
The estimate on 
each 
$\left(\Tbar\sb{\lambda,\sigma}\sp{\hvar}\right)
\sb{X\varTheta}$
is straightforward:
The mixed Hessian $S\sb{x\thvar}$ 
is of rank at least $n-1$, while the $x$-support 
of the integral kernel 
is of size $\hvar^{\frac{1}{k}}$,
and $\thvar$-support is of size $\hvar$.
Therefore, 
according to
H\"ormander's estimate 
for non-degenerate oscillatory integrals 
in $\R^{n-1}$ \cite{H\"o\ye{71}}
(in $x\pr$, $\thvar\pr$-variables)
and 
to the Schur lemma
(in $x\sb{n}$, $\thvar\sb{n}$-variables), 
we conclude that 
$$
\norm{\left(\Tbar\sb{\lambda,\sigma}\sp{\hvar}
\right)\sb{X\varTheta}}
\le\const\lambda^{-\frac{n-1}{2}}
(\hvar^{\frac{1}{k}}\hvar)^{\frac{1}{2}}.
\tag{\sec3.7}
$$
This agrees with $(\sec3.1)$.

The almost orthogonality of the pieces
localized near different points in the $\thvar$-space
is easy to establish. This orthogonality is proved
identically to the almost orthogonality of
pieces $\left(T\sb{\lambda}\sp{\hvar}\right)\sb{\varTheta}$
from Section \sec2;
again, we refer to \cite{Cu\ye{97}}.
Therefore, we can assume that
$\varTheta$ is the same
for all 
$\left(\Tbar\sb{\lambda,\sigma}\sp{\hvar}
\right)\sb{X\varTheta}$,
and we only need to prove 
the almost orthogonality with 
respect to different values of $X$.
We put for brevity
$$
\bar{\tau}\sp{}\sb{X}
\equiv\left(\Tbar\sb{\lambda,\sigma}\sp{\hvar}
\right)\sb{X\varTheta},
\qquad
\bar{\tau}\sp{}\sb{Y}
\equiv\left(\Tbar\sb{\lambda,\sigma}\sp{\hvar}
\right)\sb{Y\varTheta};
$$
the values of $\sigma$ and $\varTheta$
are assumed to be the same for the rest
of the section.

We claim that 
these operators
are almost orthogonal:

\proclaim{Proposition \sec3.2}
The operators 
$
\bar{\tau}\sp{}\sb{X}=
\left(\Tbar\sb{\lambda,\sigma}\sp{\hvar}
\right)\sb{X\varTheta}
$
are almost orthogonal
with respect to different values of $X\in\Z^n$:
$$
\norm{\bar{\tau}^{\ast}\sb{X}\bar{\tau}^{}\sb{Y}}\sb{L^2\to L^2},
\ \norm{\bar{\tau}^{}\sb{X}\bar{\tau}^\ast\sb{Y}}\sb{L^2\to L^2}
\,\le\,\const \tau^2\Abs{X-Y}^{-N},
\qquad\text{for any \ $N>0$}.
$$
Here $\tau$ is the estimate $(\sec3.7)$
which is valid for each 
operator
$\bar{\tau}\sp{}\sb{X}$.

\endproclaim

Now the statement of Theorem \sec3.1 
would follow from the Cotlar-Stein lemma.

\proof
It suffices to consider the almost orthogonality 
for $\abs{X-Y}\ge 2\sqrt{n}+1$,
when the integral kernels of 
$\bar{\tau}^{}\sb{X}$ and $\bar{\tau}^\ast\sb{Y}$
have no common support in $x$.
The almost orthogonality is straightforward
for the compositions 
$
\bar{\tau}\sp{\ast}\sb{X}\bar{\tau}\sp{}\sb{Y};
$
this leaves us with
$
\bar{\tau}\sp{}\sb{X}
\bar{\tau}\sb{Y}\sp{\ast}.
$ 

The integral kernel of 
$\bar{\tau}^{}\sb{X}\bar{\tau}^\ast\sb{Y}$ 
is given by
$$
\aligned
K(\bar{\tau}^{}\sb{X}\bar{\tau}^\ast\sb{Y})(x,y)
=\int d\thvar\,e^{i\lambda(S(x,\thvar)-S(y,\thvar))}
\times\dots\,.
\endaligned
\tag{\sec3.8}
$$

It is convenient 
to fix $\thvar$
and
to work in the space $(\eta\pr,x\sb{n})$,
which is the image of the diffeomorphism
$
\pi\sp{\prime}:\;
x\mapsto 
(\eta\pr(x)\equiv S\sb{\thvar\cpr}(x,\thvar),x\sb{n}),
$
which already appeared in $(\sec2.13)$.
We denote by  ${\Cal L}$
the line segment
from $(\eta\pr(x),x\sb{n})$ to $(\eta\pr(y),y_n)$;
$\abs{{\Cal L}}$ stays for the length of ${\Cal L}$.
Without the loss of generality we assume
$\abs{{\Cal L}}\le 1$.

We will consider two cases:

\noindent
$\bullet$
The {\it vertical case},
when the line segment ${\Cal L}$ 
is within
the conic neighborhood of magnitude
$$
\alpha=c\hvar^{1-\frac{1}{k}},
\qquad
\text{for some small $c>0$},
$$ 
of the directions 
$\pm\left(\p\sb{x\csb{n}}\right)\sb{\eta\cpr}$.

Let $t$ be a parameter on the line segment
${\Cal L}$, 
which changes
from $t=0$ at $\pi\ssb{R}^{\prime}\at{\thvar}(x)$ 
to $t=\abs{{\Cal L}}$ 
at $\pi\ssb{R}^{\prime}\at{\thvar}(y)$.
Since $\pi\sp{\prime}$ is a diffeomorphism,
we may assume that
$
c_1\abs{x-y}
\le
\abs{{\Cal L}}
\le c_2 \abs{x-y},
$
for some constants $c_2>c_1>0$,
and since
$\abs{x-y}\approx 
\hvar^{\frac{1}{k}}\Abs{Y-X}
$
(with the error of magnitude $\hvar^{\frac{1}{k}}$),
we have 
$$
C_1\hvar^{\frac{1}{k}}\abs{X-Y}
\le \abs{{\Cal L}}\le
C_2\hvar^{\frac{1}{k}}\abs{X-Y},
\qquad C_2>C_1>0.
\tag{\sec3.9}
$$

We may consider $h\at{\Cal L}$ as a function of $t$.
Since ${\Cal L}$ is in the $\alpha$-cone 
around $\pm\left(\p\sb{x\csb{n}}\right)\sb{\eta\cpr}$,
$$
\p\sb{t}^j h(t)
=\vect{K}\ssb{R}^j h(x,\thvar)\at{\Cal L}
\quad\text{modulo terms of magnitude }\, 
\alpha=c\hvar^{1-\frac{1}{k}}.
\tag{\sec3.10}
$$
Therefore,
the values of the derivatives 
$h\sp{(j)}(t)$
are close 
to the values of $\vect{K}\ssb{R}^j h(x,\thvar)\at{\Cal L}$.
Since 
$\sigma\sb{j}\vect{K}\ssb{R}\sp{j}h\ge-\hvar$
at the points 
$(x,\thvar)$ and $(y,\thvar)$
(which are on the support of the integral kernel 
of $T\sb{\lambda,\sigma}\sp{\hvar}$),
we can take $c$ small enough
so that 
for any $j<k$
$$
\sigma\sb{j}h\sp{(j)}(0)
\ge -\hvar^{1-\frac{1}{k}},
\qquad
\sigma\sb{j}h\sp{(j)}(\abs{{\Cal L}})
\ge -\hvar^{1-\frac{1}{k}}.
\tag{\sec3.11}
$$
According to $(\sec2.10)$, 
$\abs{\vect{K}\ssb{R}^k h}\ge \kappa>0$;
hence
(if $c$ is sufficiently small)
we also know that
$$
\abs{h\sp{(k)}(t)}
\ge \frac{\kappa}{2}>0,
\qquad\text{for any } t \text{ between }
0 \text{ and }\abs{{\Cal L}}.
\tag{\sec3.12}
$$
Let us show what restriction this imposes on $\abs{{\Cal L}}$.

\proclaim{Lemma}
Let $f(t)\in C^{k}(\R)$.
Assume that for some $l$,
$0\le l\le 1$, 
for some set of $k-1$ numbers
$\sigma_j=\pm 1$, \ $1\le j\le k-1$,
and for some $\epsilon>0$
the followig conditions are satisfied:
$$
\sigma\sb{j}f\sp{(j)}(0)
\ge -\epsilon
\quad\text{and}
\quad 
\sigma\sb{j}f\sp{(j)}(l)
\ge -\epsilon,
\qquad 1\le j< k,
\tag{\sec3.13}
$$
$$
\abs{f\sp{(k)}(t)}
\ge \varkappa>0
\qquad\text{for \  $0\le t\le l$}.
\tag{\sec3.14}
$$
Then 
$$
\abs{f(l)-f(0)}\ge\varkappa\frac{l^k}{k!}
-(k-1)\epsilon l.
\tag{\sec3.15}
$$

\endproclaim

Similar inequalities
appeared in \cite{Ch\ye{85}}
and \cite{PhSt\ye{97}}.

\proof
Due to $(\sec3.14)$, 
the function $f\sp{(k-1)}(t)$ is monotone.
  From $(\sec3.13)$ we know that
$\sigma\sb{k-1}f\sp{(k-1)}\ge-\epsilon$ 
at $t=0$ and $t=l$,
and, since $\abs{f\sp{(k)}(t)}\ge\varkappa$,
we derive that
$$
\text{either}\quad
\sigma\sb{k-1}
f\sp{(k-1)}(t)\ge\varkappa t
-\epsilon
\quad
\text{or}
\quad
\sigma\sb{k-1}f\sp{(k-1)(t)}\ge\varkappa(l-t)
-\epsilon,
\tag{\sec3.16}
$$ 
for any $t$ between $0$ and $l$,
depending on the relation between the signs of 
$f\sp{(k)}(t)$ and $\sigma\sb{k-1}$.

Assume that $\sigma\sb{k-1}=1$ 
and that the first inequality in $(\sec3.16)$ is 
satisfied.
If $\sigma\sb{k-2}=1$, then  $f^{(k-2)}(0)\ge-\epsilon$,
and we have:
$$
f\sp{(k-2)}(t)\ge\varkappa \frac{t^2}{2}-\epsilon t
+f\sp{(k-2)}(0)
\ge\varkappa \frac{t^2}{2}
-2\epsilon,
\qquad 0\le t\le l.
$$
If instead $\sigma\sb{k-2}=-1$, then from
$f^{(k-2)}(l)\le\epsilon$
and
$f^{(k-1)}(t)\ge\varkappa t-\epsilon$
we derive
$$
f\sp{(k-2)}(t)
\le
-\varkappa \frac{l^2-t^2}{2}+\epsilon (l-t) 
+f\sp{(k-2)}(l)
\le
-\varkappa \frac{(l-t)^2}{2}
+2\epsilon,
\qquad 0\le t\le l.
$$
All other cases are treated similarly;
each time we end up with
one of
the following bounds on $f\sp{(k-2)}(t)$:
$$
\text{either}\quad
\sigma\sb{k-2}f\sp{(k-2)}(t)\ge\varkappa \frac{t^2}{2}
-2\epsilon
\quad\text{or}\quad
\sigma\sb{k-2}f\sp{(k-2)}(t)\ge\varkappa\frac{(l-t)^2}{2}
-2\epsilon,
$$ 
depending on the relation 
between signs of $\sigma\sb{k-1}$
and $\sigma\sb{k-2}$, and which of 
the inequalities in $(\sec3.16)$
is valid.
We continue by induction 
and conclude that
$$
\text{either}\quad
\sigma\sb{1}f\pr(t)
\ge\varkappa\frac{t^{k-1}}{(k-1)!}
-(k-1)\epsilon
\quad
\text{or}
\quad
\sigma\sb{1}f\pr(t)
\ge\varkappa\frac{(l-t)^{k-1}}{(k-1)!}
-(k-1)\epsilon.
$$ 
In either case,
$
\ \dsize
\abs{f(l)-f(0)}
\ge\varkappa\frac{l^k}{k!}
-
(k-1) \epsilon l.
$
\endproof

According 
to $(\sec3.11)$, $(\sec3.12)$,
and to the above lemma 
(with $\epsilon=\hvar^{1-\frac{1}{k}}$),
$$
\abs{h(l)-h(0)}
\ge \frac{\kappa}{2}\frac{\abs{{\Cal L}}^k}{k!}
-(k-1)\hvar^{1-\frac{1}{k}}\abs{{\Cal L}}.
$$
Since the left-hand side could not be greater than $4\hvar$, 
there is the following restriction on 
the length of ${\Cal L}$:
$$
\abs{{\Cal L}}\le\const\hvar^{\frac{1}{k}}.
\tag{\sec3.17}
$$
Therefore, according to $(\sec3.9)$,
$\Abs{Y-X}\le\const$.
We conclude that 
for sufficiently large values of 
$\abs{X-Y}$
the line segment ${\Cal L}$ 
is only allowed to be outside the conic neighborhood 
of magnitude 
$c\hvar^{1-\frac{1}{k}}$
(for certain small constant $c$)
of the directions 
$\pm\left(\p\sb{x\csb{n}}\right)\sb{\eta\cpr}$
in the $(\eta\pr,x\sb{n})$-space.

\noindent
$\bullet$ 
We are thus left to consider 
the {\it the horizontal case}, when
the line segment ${\Cal L}$ is outside the 
$\alpha$-cone around 
$\pm \left(\p\sb{x\csb{n}}\right)\sb{\eta\pr}$,
where
$\alpha= c\hvar^{1-\frac{1}{k}}$
and $c>0$ is some small constant:
$$
\abs{\eta\pr(y)-\eta\pr(x)}\ge\sin\alpha
\cdot
\left(
 \abs{\eta\pr(y)-\eta\pr(x)}^2+\abs{y_n-x\sb{n}}^2
\right)^{1/2}.
\tag{\sec3.18}
$$
We use $(\sec3.9)$ and obtain 
$$
\abs{\eta\pr(y)-\eta\pr(x)}
\ge\const\alpha\hvar^{\frac{1}{k}}\Abs{Y-X}
\ge\const\hvar\Abs{Y-X}.
\tag{\sec3.19}
$$

We integrate in $(\sec3.8)$ by parts, with the aid of
the operator 
$$
L\sb{\thvar}=
\frac{1}{i\lambda}\frac{
\left(S\sb{\thvar}(x,\thvar)
-S\sb{\thvar}(y,\thvar)\right)
\cdot\habla{\thvar}}
{\Abs{S\sb{\thvar}(x,\thvar)
-S\sb{\thvar}(y,\thvar)}^2}.
$$
Each derivative $\habla{\thvar}$
contributes at most
$\hvar^{-1}$.
According to $(\sec3.19)$,
this also includes the case when the
derivative falls on the denominator of $L\sb{\thvar}$
itself.
Therefore, each integration by parts
yields the factor
$$
\frac{\const}{\lambda\hvar\cdot\hvar
\Abs{Y-X}}.
\tag{\sec3.20}
$$
According to $(\sec3.2)$, 
$\hvar
\ge\lambda^{-\frac{k}{2k+1}}$,
and therefore
$
\lambda\hvar^2\ge\lambda^{1-\frac{2k}{2k+1}}
=\lambda^{\frac{1}{2k+1}}.
$
Repeated integration by parts in $(\sec3.8)$
yields powers of $(\sec3.20)$, and we gain 
arbitrarily large negative powers of $\lambda$
and $\Abs{Y-X}$.
This proves the required almost orthogonality 
relations
for the operators $\bar{\tau}^{}\sb{X}$ with different indices $X\in\Z^n$,
and concludes the proof of Proposition \sec3.2.
\endproof

\bigskip\head
\sec4. Sharpness of the results
\endhead

Let us consider 
a particular oscillatory integral operator 
$T\sp{(1,2)}\sb{\lambda}$,
$$
T\sp{(1,2)}\sb{\lambda}\,u(x)
=\underset{\R}\to\int
e^{i\lambda S(x,\thvar)}
\psi(x,\thvar)
u(\thvar)\,d\thvar,
\qquad x,\,\thvar\in\R,
\tag{\sec4.1}
$$
with the phase function given by
$$
S(x,\thvar)=x^3\thvar-x\thvar^2.
\tag{\sec4.2}
$$
The function $\psi\in C\comp^\infty(\R\times\R)$
is supported in the unit ball centered in the origin
in $\R\times\R$.
We assume that
near the origin
$\psi\equiv 1$.

The projections from the associated canonical relation
are represented by the maps
$$
\pi\ssb{L}:\;
(x,\thvar)
\mapsto (x,S\sb{x}=3x^2\thvar-\thvar^2),
\qquad
\pi\ssb{R}:\;
(x,\thvar)
\mapsto (\thvar,S\sb{\thvar}=x^3-2x\thvar),
$$
which have the singularities 
of the Whitney fold ($k=1$) 
and the simple cusp ($k=2$),
respectively.
This is represented by the superscript
$(1,2)$.
(Note that the determinants of the Jacobi matrices
of both projections are equal to 
$h(x,\thvar)=3x^2-2\thvar$,
so that 
$\p\sb{\thvar} h\ne 0$,
$\p\sb{x}^2 h\ne 0$.)
According to Theorem \sec1.2,
$
\norm{T\sp{(1,2)}\sb{\lambda}}
\le\const\lambda^{-\frac{1}{2}+\frac{1}{5}}.
$
We are going to prove that this estimate is optimal.

\proclaim{Proposition \sec4.1}
The optimal rate of decay
of $\norm{T\sp{(1,2)}\sb{\lambda}}\sb{L^2\to L^2}$ 
equals $3/10$.
\endproclaim

Let us assume that the operator 
$T\sp{(1,2)}\sb{\lambda}$ is bounded
from $L^2$ to $L^2$ by
$$
\norm{T\sp{(1,2)}\sb{\lambda}}\le\const\lambda^{-d},
\tag{\sec4.3}
$$
where $d$ is some positive real number
(which {\it a priori} could be greater
than $3/10$).
We consider the family of operators,
$$
T\sp{(1,2)}\sb{\lambda,R}\,u(x)
=\underset{\R}\to\int
e^{i\lambda S(x,\thvar)}
\psi\left(\frac{x}{R},\frac{\thvar}{R^2}\right)
u(\thvar)\,d\thvar,
\qquad x,\,\thvar\in\R,
\tag{\sec4.4}
$$
where $R\ge 1$.
All these operators 
are bounded from $L^2$ to $L^2$
(as long as $R<\infty$).
We would like to know the behavior of
their norms as $\lambda$ and $R$ become large.

We rescale $x$ and $\thvar$
with the aid of some $\mu>0$:
$$
\aligned
T\sp{(1,2)}\sb{\lambda,R}\,u(\mu x)
&=\underset{\R}\to\int
e^{i\lambda S(\mu x,\mu^2 \thvar)}
\psi\left(\frac{\mu x}{R},\frac{\mu^2 \thvar}{R^2}
\right)
u(\mu^2\thvar)\,d(\mu^2\thvar)
\\
&=\mu^2\underset{\R}\to\int
e^{i\lambda\mu^5 S(x,\thvar)}
\psi\left(\frac{\mu x}{R},\frac{\mu^2 \thvar}{R^2}
\right)
u(\mu^2\thvar)\,d\thvar.
\endaligned
$$
We put $\mu=R$
and use the assumption $(\sec4.3)$,
getting
$$
\norm{T\sp{(1,2)}\sb{\lambda,R}u(R x)}\sb{L^2}
\le \const R^2(\lambda R^5)^{-d}
\norm{u(R^2 \thvar)}\sb{L^2}.
\tag{\sec4.5}
$$

Now we rescale $x$ and $\thvar$ ``back''
and keep track of the powers of $R$,
obtaining
$$
R^{-\half}\norm{T\sp{(1,2)}\sb{\lambda,R}u(x)}\sb{L^2}
\le\const R^2(\lambda R^5)^{-d} 
R^{-1}\norm{u(\thvar)}\sb{L^2},
$$
which gives the following bound on
$T\sp{(1,2)}\sb{\lambda,R}$:
$$
\norm{T\sp{(1,2)}\sb{\lambda,R}}\le 
R^{\frac{3}{2}-5d}\const\lambda^{-d}.
\tag{\sec4.6}
$$

Now let us argue 
that the exponent $d=3/10$ is optimal.
Assuming $d>3/10$,
we could conclude from $(\sec4.6)$ that
$
\norm{T\sp{(1,2)}\sb{\lambda,R}}\to 0
$
as $R$ becomes large (and $\lambda$ is fixed).
At the same time,
if we take a function $u(\thvar)$ supported 
in a small neighborhood of $\thvar=0$,
then 
the image $T\sp{(1,2)}\sb{\lambda,R}u(x)$
would not change the values, 
in some small neighborhood of $x=0$,
when $R$ grows up.
Therefore,
$\norm{T\sp{(1,2)}\sb{\lambda,R}u(x)}\sb{L^2}$
could not decrease,
and we are facing the contradiction.

\demo{Remark 1}
A slight modification of the proof shows that
the decay $\sim\lambda^{-0.3}$ 
is the sharp result for the decrease 
of the norm of $T\sp{(1,2)}\sb{\lambda}$,
in the sense that 
for each $\lambda$
one can choose a function 
$u\sp{(\lambda)}\in C\comp^\infty(\R)$
supported in a small neighborhood of the origin 
such that
$$
\norm{T\sp{(1,2)}\sb{\lambda}u\sp{(\lambda)}(x)}\sb{L^2}
\ge c\lambda^{-0.3}\norm{u\sp{(\lambda)}(x)}\sb{L^2},
$$
with the constant $c>0$ independent on $\lambda$.
\enddemo

\demo{Remark 2}
Since $(\sec4.6)$ 
is valid with $d=3/10$,
the norm of the operator
$T\sp{(1,2)}\sb{\lambda,R}$
defined by $(\sec4.4)$
does not increase as $R$ becomes large:
$
\norm{T\sp{(1,2)}\sb{\lambda,R}}
\le\const\lambda^{-3/10}
$
independent on a particular value of $R$.
Hence, the operators
$T\sp{(1,2)}\sb{\lambda,R}$ converge
(in the weak $L^2\to L^2$
operator topology)
to the non-compactly supported 
oscillatory integral operator
$\tilde{T}\sp{(1,2)}\sb{\lambda}$,
defined by
$$
\tilde{T}\sp{(1,2)}\sb{\lambda}\,u(x)
=\underset{\R}\to\int
e^{i\lambda (x^3\thvar-x\thvar^2)}
u(\thvar)\,d\thvar,
\qquad x,\,\thvar\in\R.
\tag{\sec4.7}
$$
$\tilde{T}\sp{(1,2)}\sb{\lambda}$
extends to a continuous operator 
on $L^2$
with the norm
$\const\lambda^{-3/10}$.
\enddemo

\subhead{Models of 
operators
with higher order singularities}
\endsubhead

We generalize the previous example
and construct the canonical relation
with the projection $\pi\ssb{L}$ 
being a fold
and $\pi\ssb{R}$ being a map 
with a Morin $\sing{1\sb{k}}$-singularity.

We fix $n\ge k-1$
and introduce the phase function
$S(x,\thvar)\in C^\infty(\R^n\times\R^n)$
given by the polynomial
$$
\aligned
S(x,\thvar)
&=(x\sb{n}^{k+1}+x\sb{n}^{k-1}x\sb{n-1}
+\dots
+x\sb{n}^2 x\sb{n-k+2})\thvar_n
+x\sb{n}\frac{\thvar_n^2}{2}
+x\pr\cdot\thvar\pr,
\\
&
\qquad x=(x\pr,x\sb{n})\in\R^n,
\quad
\thvar=(\thvar\pr,\thvar_n)\in\R^n.
\endaligned
\tag{\sec4.8}
$$

The map $\pi\ssb{R}$, 
$$
\left[\matrix
x\pr
\\
x\sb{n}
\\ \\
\thvar\pr \\
\thvar_n
\endmatrix\right]
\overset{\pi\ssb{R}}\to\longmapsto
\left[\matrix
\thvar\pr
\\
\thvar_n
\\ \\
S\sb{\thvar\cpr} \\
S\sb{\thvar\csb{n}}
\endmatrix\right]
=
\left[\matrix
\thvar\pr\hfill
\\ 
\thvar_n\hfill\\ 
\\
x\pr \hfill\\ 
x\sb{n}^{k+1}+x\sb{n}^{k-1}x\sb{n-1}+\dots+x\sb{n}^2 x\sb{n-k+2}
+x\sb{n}\thvar_n
\endmatrix\right],
$$
has the canonical form \cite{Mo\ye{65}}
of a map with a Morin 
$\sing{1\sb{k}}$-singularity 
at the origin.
Then, since 
$\det d\pi\ssb{L}=\det d\pi\ssb{R}
=(k+1)x\sb{n}^k+\dots+\thvar_n$
vanishes of the first order in the direction 
of the kernel of $\pi\ssb{L}$ 
(which is generated 
at $x=\thvar=0$
by $\p\sb{\thvar\csb{n}}$),
$\pi\ssb{L}$ is a Whitney fold.

We consider the oscillatory integral 
operator
with the phase function $(\sec4.8)$,
$$
T\sp{(1,k)}\sb{\lambda}u(x)
=\underset{\R^n}\to\int
e^{i\lambda S(x,\thvar)}
\psi(x,\thvar)
u(\thvar)\,d\thvar,
\qquad x,\,\thvar\in\R^n,
\tag{\sec4.9}
$$
where $\psi$ is a smooth function
supported near the origin.
According to Theorem \sec1.2,
$$
\norm{T\sp{(1,k)}\sb{\lambda}}
\le\const\lambda^{-\frac{n}{2}+(4+\frac{2}{k})^{-1}}.
$$

\proclaim{Proposition \sec4.2}
The optimal rate of decay of 
$\norm{T\sp{(1,k)}\sb{\lambda}}$
equals 
$\,\frac{n}{2}-(4+\frac{2}{k})^{-1}$.
\endproclaim

\proof
The proof is similar to the proof of 
Proposition \sec4.1.
If we rescale $x\sb{n}\mapsto \mu x\sb{n}$,
then 
for $S$ to be homogeneous in $\mu$
we need to rescale $x$ and $\thvar$
as follows:
$$
x\to X\sb{\mu}(x)
=(
\mu^{n}x_1,\,
\dots,\,
\mu^{n-j+1} x\sb{j},\,
\dots,\,
\mu x\sb{n}),
$$
$$
\thvar\to \Theta\sb{\mu}(\thvar)
=(
\mu^{2k+1-n}\thvar_1,\,
\dots,\,
\mu^{2k+j-n}\thvar_j,\,
\dots,\,
\mu^{2k-1}\thvar\sb{n-1},\,
\mu^k \thvar_n).
$$
We then have
$
\dsize
S(X\sb{\mu}(x),\Theta\sb{\mu}(\thvar))
=\mu^{2k+1}S(x,\thvar).
$
We define
$$
T\sb{\lambda,R}\sp{(1,k)}
u(x)=\int\sb{\R^n}
e^{i\lambda S(x,\thvar)}
\psi(X\sb{R^{-1}}(x), \Theta\sb{R^{-1}}(\thvar))
u(\thvar)d\thvar.
$$
%
We proceed similarly to the proof of Proposition \sec4.1
and obtain
$$
\norm{T\sb{\lambda,R}\sp{(1,k)}}
\le\const (\lambda R^{2k+1})^{-d}
\Abs{\frac{\p X\sb{R}}{\p x}}^{\half}
\Abs{\frac{\p \Theta\sb{R}}{\p \thvar}}^{\half}.
$$
Here $\Abs{\frac{\p X\sb{R}}{\p x}}$,
$\Abs{\frac{\p \Theta\sb{R}}{\p\thvar}}$
are the determinants of the Jacobi matrices
of the maps
$X\sb{\mu}(x)$ and $\Theta\sb{R}(\thvar)$
(which only depend on $R$).
To simplify the rest, we notice that
$$
\frac{\p X\sb{R}}{\p x}\cdot
\frac{\p \Theta\sb{R}}{\p \thvar}
=\operatorname{diag}
(R^{2k+1},\,\dots,\,R^{2k+1},\,R^{k+1}),
$$
$$
\Abs{
\frac{\p X\sb{R}}{\p x}\cdot
\frac{\p \Theta\sb{R}}{\p \thvar}
}=R^{n(2k+1)-k},
$$
and hence
$$
\norm{T\sb{\lambda,R}\sp{(1,k)}}
\le\const (\lambda R^{2k+1})^{-d}
R^{\frac{n(2k+1)-k}{2}}.
$$
Since the norm
of 
$T\sb{\lambda,R}\sp{(1,k)}$ 
can not decrease when $R$
becomes large
(according to the same arguments 
as in the proof of Proposition \sec4.1),
we conclude that
the rate of decay $d$
can not be larger than
$\frac{n}{2}-\frac{k}{2(2k+1)}$.
\endproof

\bigskip\head
ACKNOWLEDGMENTS
\endhead

The author is grateful to 
S. Cuccagna, A. Greenleaf, D.H. Phong,
and M. Zworski 
for their interest 
and for fruitful discussions.


\Refs\widestnumber\key{AAAA\ye{91}}

\noref\key{ArGuVa\ye{85}}
\by V.I. Arnold, S.M. Gusein-Zade, A.N. Varchenko
\book Singularities of differentiable maps. Vol. I
\publ Birkhauser Boston, Inc., Boston, MA \yr 1985
\endref

\noref\key{BaSy\ye{91}}
\by G. Bao, W.W. Symes
\paper A trace theorem 
for solutions of linear partial differential
equations
\jour  Math. Methods Appl. Sci. \vol 14 \yr 1991
\pages 553--562
\endref

\ref\key{Ch\ye{85}}
\by M. Christ
\paper 
Hilbert transforms along curves
\jour Ann. of Math. \vol 122 \yr 1985
\pages 575--596
\endref

\noref\key{Ch\ye{95}}
\by M. Christ
\paper 
Failure of an endpoint estimate for integrals along curves
\inbook 
Fourier analysis
and partial differential equations 
(Miraflores de la Sierra, 1992) \pages 163--168
\publ Stud. Adv. Math., CRC 
\publaddr Boca Raton, FL \yr 1995
\endref

\noref\key{Co\ye{97}}
\by A. Comech 
\paper 
Oscillatory integral operators
in Scattering Theory
\jour Comm. Partial Differential Equations
\vol 22 \yr 1997 \pages 841--867
\endref

\ref\key{Co\ye{97}}
\by A. Comech
\paper
Integral operators
with singular canonical relations
\inbook
Spectral theory, microlocal analysis,
singular manifolds
\eds M. Demuth, E. Schrohe, 
B.-W. Schulze, and J. Sj\"ostrand
\publ Akademie Verlag \publaddr Berlin \yr 1997
\pages 200--248
\endref

\noref\key{Co\ye{97}}
\bysame
\paper
Optimal estimates
for Fourier integral operators
with one-sided folds
\paperinfo preprint
\yr 1997
\endref

\noref\key{Co\ye{97}}
\by A. Comech
\paper
Damping estimates
for
oscillatory integral operators
with finite type singularities
\paperinfo preprint,
http://www.fields.toronto.edu/\~{}acomech/papers/chi.dvi
\yr 1997
\endref

\noref\key{Co\ye{98}}
\bysame
\paper
Sobolev estimates 
for Radon transform of Melrose and Taylor
\jour Comm. Pure Appl. Math \vol 51 \yr 1998 \pages 537--550 
\endref

\ref\key{Cu\ye{97}}\by S. Cuccagna
\paper 
$L^2$ estimates for  averaging operators  
along curves with two sided $k$ fold singularities
\jour Duke Journal
\vol 89 \yr 1997 \pages 203--216
\endref

\noref\key{F}\by B.V. Fedosov
\book Deformation quantization and index theory
\publ Akademie Verlag, Berlin \yr 1996
\endref

\noref\key{De\ye{97}}
\by J.M. Delort
\paper Sur le temps d'existence 
pour l'equation de Klein-Gordon semi-lineaire
en dimension 1
\jour Bull. Soc. math. France \vol 125 \yr 1997
\pages 269--311
\endref

\noref\key{Fw}\by D. Fujiwara
\paper On the boundedness of integral transformations 
with highly oscillatory kernels
\jour Proc. Japan Acad. \vol 51 \yr 1975 \pages 96--99
\endref

\ref\key{GrSe\ye{94}}\by A. Greenleaf and A. Seeger
\paper  Fourier integral operators with fold singularities
\jour J. Reine Angew. Math. \vol 455 \yr 1994 
\pages 35--56
\endref

\ref\key{GrSe\ye{98}}\bysame
\paper  Fourier integral operators with cusp singularities
\jour Amer. J. Math. \vol 120 \yr 1998 \pages 1077-1119 
\endref

\noref\key{GrSe\ye{97}}\by A. Greenleaf and A. Seeger
\paper  On oscillatory integral operators
with folding canonical relations
\paperinfo preprint
\yr 1997
\endref

\noref\key{GrU\ye{91}}
\by A. Greenleaf and G. Uhlmann
\paper
Composition of some singular Fourier 
integral operators and estimates for 
restricted X-ray transforms. II
\jour  Duke Math. J. \vol  64 \yr 1991 \pages 415--444
\endref

\ref\key{H\"o\ye{71}}\by L. H\"ormander 
\paper Fourier integral operators 
\vol 127 \yr1971 \pages 79--183
\jour Acta Math.\endref

\ref\key{H\"o\ye{85}}\bysame
\book The analysis of linear 
partial differential operators
\publ Springer-Verlag \yr 1985 
\endref

\noref\key{KlMa\ye{96}}
\by S. Klainerman, M. Machedon
\paper Remark on Strichartz-type inequalities
\jour Internat. Math. Res. Notices \yr 1996\pages 201--220
\endref

\noref\key{KoVa\ye{96}}
\by A.I. Komech, B.R.  Vainberg
\paper On asymptotic stability 
of stationary solutions to nonlinear 
wave and Klein-Gordon equations
\jour Arch. Rational Mech. Anal. \vol 134 \yr 1996
\pages 227--248
\endref

\ref\key{Me\ye{76}}
\by R.B. Melrose
\paper Equivalence of glancing hypersurfaces
\jour Invent. Math. \vol 37 \yr 1976
\pages 165--191
\endref

\noref\key{MeSj\ye{78}}
\by R.B. Melrose and J. Sj\"ostrand
\paper Singularities of boundary value problems, I
\jour Comm. Pure Appl. Math. \vol 31 \yr 1978
\pages 593--617
\endref

\ref\key{MeTa\ye{85}}\by R.B. Melrose and M.E. Taylor 
\paper Near peak scattering and the corrected
Kirchhoff approximation for a convex obstacle
\jour Adv. in Math. \vol 55 \yr 1985 \pages 242--315
\endref 

\ref\key{Mo\ye{65}}
\by B. Morin
\paper 
Canonical forms of the singularities of a differentiable mapping
\jour C. R. Acad. Sci. Paris \vol 260 \yr 1965 
\pages 6503--6506
\endref

\noref\key{Pa\ye{91}}\by Y.B. Pan 
\paper Uniform estimates for oscillatory integral operators 
\vol 100\yr 1991\pages 207--220
\jour J. Funct. Anal. \endref

\ref\key{PaSo\ye{90}}\by Y.B. Pan and C.D. Sogge
\paper Oscillatory integrals associated to folding 
canonical relations
\jour Colloq. Math. \vol 61 \yr1990 \pages 413--419 
\endref

\ref\key{Ph\ye{94}}\by D.H. \!Phong
\paper Singular integrals and Fourier integral operators
\inbook Essays on Fourier Analysis in honor of Elias M. Stein
\eds C. Fefferman, R. Fefferman and S. Wainger
\publ Princeton Univ. Press\yr 1994 \pages 287--320 
\endref

\noref\key{Ph\ye{95}}\bysame
\paper Regularity of Fourier integral operators
\jour Proceedings of the International Congress
of Mathematicians\vol 1, 2\yr 1994 \pages 862--874
\endref

\noref\key{PhSt\ye{86}}\by D.H. Phong and E.M. Stein 
\paper Hilbert integrals, Singular integrals, 
and Radon transforms I
\jour Acta Math. \vol 157 \yr1986 \pages 99--157 \endref

\ref\key{PhSt\ye{91}}\by D.H. Phong and E.M. Stein 
\paper Radon transform and torsion 
\jour Internat. Math. Res. Notices\vol 4\yr 1991\pages 49--60
\endref

\noref\key{PhSt\ye{92}}\bysame
\paper  Oscillatory integrals with polynomial phases
\jour  Invent. Math. \vol 110 \yr 1992 \pages 39-62\endref

\noref\key{PhSt}\by D.H. Phong and E.M. Stein 
\paper On a stopping process for oscillatory integrals 
\jour J. Geometric Analysis \vol 4 \yr 1994 \pages 105--120
\endref

\noref\key{PhSt}\by D.H. Phong and E.M. Stein 
\paper Operator versions of the Van der Corput lemma 
and Fourier integral operators
\jour Math. Res. Letters \vol 1 \yr 1994 \pages 27--33
\endref

\ref\key{PhSt\ye{94}}\bysame
\paper Models of degenerate Fourier integral operators 
and Radon transforms 
\jour Ann. of Math. \vol 140 \yr 1994 \pages 703--722
\endref

\ref\key{PhSt\ye{97}}\bysame
\paper Newton polyhedron and oscillatory integral operators
\jour Acta Math \vol 179 \yr 1997 \pages 105--152
\endref

\ref\key{Se\ye{93}}\by A. Seeger
\paper Degenerate Fourier integral operators in the plane
\jour Duke Math. J. \vol 71 \yr 1993\pages 685--745\endref

\ref\key{Se\ye{98}}\bysame
\paper Radon transforms and finite type conditions
\jour Journal Amer. Math. Soc. \vol 11 \yr 1998 \pages 869-897 
\endref

\noref\key{SeSoSt\ye{91}}
\by A. Seeger, C.D. Sogge and E.M. Stein
\paper  Regularity properties of Fourier integral
operators \jour Ann. of Math. 
\vol 134\yr 1991 \pages 231--251
\endref

\noref\key{SmSo\ye{94}}\by H.F. Smith and C.D. Sogge 
\paper $L^p$ regularity for the wave equation 
with strictly convex obstacles 
\jour Duke Math. J.\vol 73  \yr  1994 \pages 97--153
\endref

\noref\key{So\ye{93}}
\by C.D. Sogge
\book Fourier integrals in classical analysis
\publ Cambridge University Press
\publaddr Cambridge
\yr 1993
\endref

\noref\key{SoSt\ye{86}}
\by C.D. Sogge and E.M. Stein
\paper  Averages over hypersurfaces II
\jour Invent. Math. \vol 86\yr 1986 
\pages 233--242
\endref

\ref\key{St\ye{93}}\by  E.M. Stein
\book Harmonic Analysis: real-variable methods,
orthogonality, and oscillatory integrals 
\publ Prince\-ton Univ. Press \yr 1993 \endref

\noref\key{Ta\ye{96}}\by  D. Tataru
\paper 
The $X\sp{s}\sb{\theta}$ spaces
and unique continuation
for solutions
to the semilinear wave equation
\jour Comm. Partial Differential Equations
\vol 21 \yr 1996 \pages 841--887
\endref

\noref\key{Ta\ye{97}}\by  D. Tataru
\paper On the regularity of boundary traces
for the wave equation
\paperinfo preprint \yr 1997
\endref

\ref\key{Th\ye{63}}
\by R. Thom
\paper 
Les singularit\'es
des applications differentiables
\jour Ann. Inst. Fourier
\vol 6 \yr 1963 
\pages 43--87
\endref

\endRefs

\enddocument